\documentclass[10pt]{amsart}

\usepackage[english]{babel}
\usepackage[utf8]{inputenc}
\usepackage{amsmath}
\usepackage{amssymb}
\usepackage{amsfonts}
\usepackage{amsthm}
\usepackage{mathrsfs}
\usepackage[all]{xy}
\usepackage[pdftex]{graphicx}
\usepackage{color}
\usepackage{cite}
\usepackage{url}
\usepackage{indent first}
\usepackage[labelfont=bf,labelsep=period,justification=raggedright]{caption}
\usepackage[english]{babel}
\usepackage[utf8]{inputenc}
\usepackage{hyperref}
\usepackage[colorinlistoftodos]{todonotes}
\usepackage{caption}
\usepackage{subcaption}
\usepackage[capitalize,nameinlink]{cleveref}

\usepackage[margin=1in]{geometry} 
\makeatletter
\newtheorem*{rep@theorem}{\rep@title}
\newcommand{\newreptheorem}[2]{%
\newenvironment{rep#1}[1]{%
 \def\rep@title{#2 \ref{##1}}%
 \begin{rep@theorem}}%
 {\end{rep@theorem}}}
\makeatother

\newcommand{\eps}{\varepsilon}

\newcommand{\EE}{\mathbb{E}}

\newcommand{\NN}{\mathbb{N}}
\newcommand{\PP}{\mathbb{P}}

\newcommand{\RR}{\mathbb{R}}
\newcommand{\ZZ}{\mathbb{Z}}

\renewcommand{\P}{\mathbf{P}}

\makeatletter
\newcommand\thankssymb[1]{\textsuperscript{\@fnsymbol{#1}}}
\makeatother

\theoremstyle{plain}
\newtheorem{thm}{Theorem}
\newreptheorem{theorem}{Theorem}
\newtheorem{lemma}[thm]{Lemma}
\newtheorem{cor}[thm]{Corollary}

\newtheorem{prop}[thm]{Proposition}

\newtheorem{qn}[thm]{Question}

\newtheorem{obs}[thm]{Observation}

\theoremstyle{definition}
\newtheorem{defn}[thm]{Definition}

\theoremstyle{remark}
\newtheorem{rem}[thm]{Remark}

\numberwithin{equation}{section}
\numberwithin{thm}{section}

\title{Maximum gaps in one-dimensional hard-core models}

\author{Dingding Dong\thankssymb{1}}
\thanks{\thankssymb{1} Department of Mathematics, Harvard University, Cambridge, MA. Email: \href{mailto:ddong@math.harvard.edu} {\nolinkurl{ddong@math.harvard.edu}}. }
\author{Nitya Mani\thankssymb{2}}
\thanks{\thankssymb{2} Department of Mathematics, Massachusetts Institute of Technology, Cambridge, MA. Email: \href{mailto:nmani@mit.edu} {\nolinkurl{nmani@mit.edu}}}

\begin{document}
\maketitle

\begin{abstract}
We study the distribution of the maximum gap size in one-dimensional hard-core models. First, we randomly sequentially pack rods of length $2$ onto an interval of length $L$, subject to the hard-core constraint that rods do not overlap. We find that in a saturated packing, with high probability there is no gap of size $2 - o(1/L)$ between adjacent rods,  but there are gaps of size at least $2 - 1/L^{1-\eps}$ for all $\eps > 0$.

We subsequently study a variant of the hard-core process, the one-dimensional ``ghost" hard-core model introduced by Torquato and Stillinger \cite{TS06}. In this model, we randomly sequentially pack rods of length $2$ onto an interval of length $L$, such that placed rods neither overlap with previously placed rods \textit{nor} previously considered candidate rods. We find that in the infinite time limit, with high probability the maximum gap between adjacent rods is smaller than $\log L$ but at least $(\log L)^{1-\eps}$ for all $\eps > 0.$ 
\end{abstract}

\section{Introduction}
The \textit{R\'enyi parking problem} is a classical combinatorial question that gives a simple example of a \textit{random sequential addition (RSA) process}; it is a specific instantiation of a one-dimensional \textit{hard-core model} of much interest in statistical mechanics.

The setup for the parking problem proceeds as follows. Consider a closed interval $[0, L]$ for $L > 2$, into which rods of length 2 sequentially arrive at integer times. When each rod arrives, we attempt to place it uniformly at random in the interval, subject to the \textit{hard-core} condition that rods cannot overlap with each other. In 1958, R\'enyi proved the following well-known result.

\begin{thm}[R\'enyi~\cite{REN58}]
In the above setup, let $N(L)$ be the random variable representing the number of rods placed in a saturated packing of $[0, L]$ (when no more rods can fit without violating the hard-core constraint). Then,
$$\lim_{L\rightarrow\infty}\frac{2\EE[N(L)]}{L}=\alpha,$$
where $\alpha$ is the \textit{R\'enyi parking constant}
$$\alpha:=\int_0^\infty\exp\left[-2\int_0^x\frac{1-e^{-y}}{y}dy\right]\,dx \approx 0.7475979202.$$
\end{thm}

In this work, we study the distribution of \textit{gaps} between adjacent rods in the saturated state, focusing on the upper extreme. In particular, we seek to understand the following:

\begin{qn}
What can we say about the largest gap that arises in a saturated packing of length $2$ rods into an interval of length $L$ by rods of length 2, subject to the hard-core constraint? 
\end{qn}

Itoh~\cite{ITO80} studied a delay integral equation that characterizes the distribution of the minimum gap sizes in a saturated configuration, following methods of Dvoretzky and Robbins~\cite{DR64}. The distribution of gap sizes was also examined in the study of the nearest neighbors problem in one-dimensional random sequential adsorption~\cite{RTT}.
In~\cite{ITO80}, Itoh observed that the expected minimum gap size in a saturated configuration on an interval of length $L$ is smaller than any constant $\eps > 0$ in the large $L$ limit.  This work was subsequently extended to give approximations of the upper tail of the distribution of minimum gap sizes in~\cite{MOR87}. These works, as noted in \S4 of~\cite{ITOH11}, imply an analogous integral recurrence for the CDF of the maximum gap size in a saturated packing. In~\cite{ITOH11} the authors provided some preliminary observations and noted that further study of the maximum gap size was of substantial interest.

In this work, we give a threshold for the maximum gap size in a saturated configuration of the hard-core model, observing that with high probability, a saturated one-dimensional hard-core packing on an interval of length $L$ has no gap of size $2 - o(1/L)$ but does have gaps of size $2 - 1/L^{1-\eps}$ for all $\eps > 0$. More precisely, we prove the following result.

\begin{thm}\label{t:maxgap}
The following holds in the saturated configuration of a one-dimensional  hard-core process on an interval of length $L$ packed by rods of length $2$, for $L$ sufficiently large:

\begin{itemize}
    \item with high probability, there are no gaps of size $2 - o(1/L);$
    \item for all $a>0$, with positive probability, there exists a gap of size at least $2 - a/L$;
    \item for all $\eps > 0$, with high probability,  there exists a gap of size at least $2 - 1/L^{1-\eps}$.
\end{itemize}
\end{thm}

In the classical one-dimensional hard-core model described above, cars that fail to park have no effects on future parking attempts. We will be interested in a variant of the hard-core model, motivated by the \textit{ghost RSA process} introduced in work of~\cite{TS06} studying sphere packings. Unlike the classical random sequential addition process, where much is unknown even in $2$ dimensions, the authors of~\cite{TS06} are able to analytically derive the $n$-point correlation functions and limiting densities, exactly solving the ghost sphere packing model in arbitrary dimension.

We study the \textit{one-dimensional ghost hard-core model}, akin to the hard-core model above, focusing on properties of this process in the infinite time limit. We give a precise definition below:
\begin{defn}\label{d:ghost}
We attempt to place rods of length $2$ on an interval of length $L$, as follows:
\begin{itemize}
    \item Initialize $X=[0,L]$ and $Y=\emptyset$.
    \item For $t=1,2,\dots$:
    \begin{itemize}
        \item If the $X\setminus Y$ has no connected component of length 2, abort.
        \item Choose a uniformly random point $x \in [0, L]$.
        \item If $x \in [0, 1) \cup (L - 1, L]$, reject $x$, replace $Y$ by $Y\cup ((x-1,x+1)\cap [0,L])$, and continue.
        \item If $(x-1,x+1)\cap Y\neq \emptyset$, reject $x$, replace $Y$ by $Y\cup (x-1,x+1)$, and continue.
        \item Else, accept $x$, replace $Y$ by $Y\cup (x-1,x+1)$, and continue.
    \end{itemize}
\end{itemize}

\end{defn}
Several further observations about the occupancy probabilities and the pair correlation function associated to this process can be found in~\cref{s:basicghost}.

In the infinite time limit of the one-dimensional ghost hard-core process (on an interval of length $L$), with high probability there is no gap of size $\log L$, but there are gaps of size at least $(\log L)^{1- \eps}$ for arbitrarily small $\eps > 0$. More precisely, we have the following: 

\begin{thm}\label{t:ghostgap}
The following holds in the infinite time limit of a one-dimensional ghost hard-core process on an interval of length $L$ packed by rods of length $2$, for $L$ sufficiently large:

\begin{itemize}
    \item with high probabiity, all gaps are smaller than $\log L;$
    \item for all $\eps > 0,$ with high probability, there exists a gap of size $(\log L)^{1- \eps}$.
\end{itemize}
\end{thm}

\subsection*{Overview of article}
We begin in~\cref{s:prelim} by reviewing the classical one-dimensional hard-core model and introducing the ghost RSA process of~\cite{TS06}. In~\cref{s:maxgap} we prove~\cref{t:maxgap}. We prove~\cref{t:ghostgap} in~\cref{s:ghostgap}. The ghost hard-core process is very different from the classical hard-core process; we illustrate some differences to give some context in~\cref{s:basicghost}.

One of our primary motivations for studying large gaps in these hard-core processes is to provide a glimpse into what gaps might look like in a random sequential addition process in higher dimensions.~\cref{t:ghostgap} hints that in higher dimensions, in the infinite time limit, a ghost packing may still have room for many more spheres/cubes to be packed without overlap; we discuss further in~\cref{s:discuss}.

\subsection*{Acknowledgements}
We would like to thank Henry Cohn and Salvatore Torquato for helpful discussions, ideas for writing improvements, and several useful reference suggestions for background. NM was supported by the Hertz Graduate Fellowship and  by the NSF GRFP \#2141064.

\section{Preliminaries}\label{s:prelim}
\subsection{Notation}
 Throughout this article, we consider packing rods of length $2$ onto an interval of length $L$, which we model by the closed interval $[0,L] \subseteq \RR$. Unless stated otherwise, we study configurations in the \textit{infinite time limit} of the two processes, the \textit{1D classical hard-core model}  and the \textit{1D ghost hard-core model}.
 
 We sometimes refer to the infinite time limit of the 1D classical hard-core model as \textit{saturation,} since at this limit, no more rods can be packed without violating the hard-core constraint. We sometimes omit the modifier \textit{hard-core} when describing models, as all models considered in this article are subject to the hard-core constraint. We also employ the following notation conventions. 
\begin{itemize}
    \item Let $N(L)$ denote the number of rods in the classical model at saturation and $\widetilde N(L)$ denote the number of rods in the ghost model in the infinite time limit.
    \item Let $G(L,r)$ denote the number of gaps of length at least $r$  in the classical model at saturation, and $\widetilde G(L,r)$ denote the number of gaps of length at least $r$ in the ghost model in the infinite time limit.
    \item For points $x_1,\dots,x_n$ on the interval, let $\pi(x_1,\cdots,x_n;L)$ be the \textit{$n$-point correlation function} in the classical hard-core model, the probability that all of $x_1,\dots,x_n$ are occupied by rods at saturation; let $\widetilde\pi(x_1,\cdots,x_n;L)$ be the corresponding $n$-point correlation function of the ghost model in the infinite time limit.
\end{itemize}

For quantities depending on $L$, we use $f=o(g)$, $f\ll g$ and $g\gg \ell$ interchangeably to denote that $\lim_{L\to\infty}f/g=0$; we use $f=O(g)$ to denote that there exists a constant $C\geq 0$ such that $f\leq Cg$ for sufficiently large $L$; we use $f\sim g$ to denote that $\lim_{L\to\infty}f/g=1$

\subsection{Classical 1D hard-core model at saturation}\label{s:1dhc}
Consider the classical 1D hard-core model, where we place rods of length $2$ on an interval of length $L$. It is easy to check the following recurrence relation on $\EE[N(L)]$:
\footnotesize 
\begin{align*}
\EE[N(L)] &= \int_{x = 1}^{L-1} \frac{1}{L-2}(1 + \EE[N(x-1)]  + \EE[N(L-x-1)])  dx =1 + \frac{2}{L-2} \int_0^{L-2} \EE[N(x)]\,dx.
\end{align*}
\normalsize 
As noted in the introduction, in~\cite{REN58}, R\'enyi established that this mean density of rods converges to the R\'enyi parking constant $\alpha \approx 0.748$.
Dvoretzky and Robbins~\cite{DR64} gave a more refined estimate of the rate of convergence. In particular, they proved that
\[
\EE[N(L)]=\alpha L/2+\alpha-1+O((4e/L)^{L/2-3/2}),
\]
indicating very fast convergence of the expected parking density $2\EE[N(L)]/L$ to the following approximate density $\alpha+\frac{2\alpha-2}{L}$.

\begin{figure}[t]
    \centering
    \includegraphics[width=0.5\textwidth]{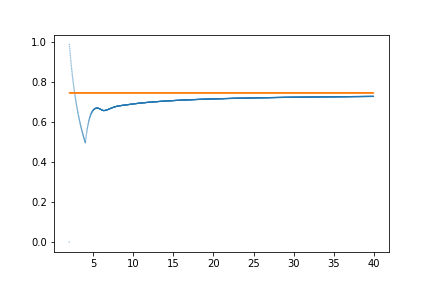}\includegraphics[width=0.5\textwidth]{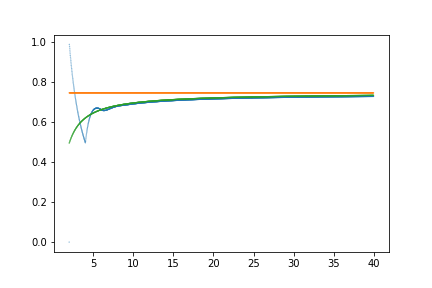}
        \caption{R\'enyi's parking constant. The blue curve denotes the expected parking density $2\EE[N(L)]/L$, the orange curve denotes the R\'enyi parking constant $\alpha$, and the green curve (depicted on the right) denotes the approximate function $\alpha+\frac{2\alpha-2}{L}$.}
    \label{fig:renyi-constant}
\end{figure}

Understanding the \textit{n}-point correlation functions is a primary motivating question when analyzing statistical mechanics models. The occupancy probability, $\pi(x,L)$ is the chance that point $x$ is covered by a rod at saturation; it has the basic symmetry property $\pi(x, L) = \pi(L-x, L)$ for all $x \in [0, L]$ and was studied along with other statistics of the correlation function in the one-dimensional hard-core model in~\cite{DH03}. 

One similar observation that arises from such analysis concerns the pair correlation function, $\pi(x_1, x_2; L)$, the probability that both $x_1, x_2$ are occupied at saturation. It will be convenient to think about this correlation when both $x_1, x_2$ are far away from the boundary, and thus we could imagine packing rods of length $2$ on a circle of length $L$. We observe that this pair correlation function is identical to the cover probability on an interval of length $L - 2$, by cutting one of the rods in half and unwinding, i.e. 
$$\pi(x_1, x_2; L) = \pi(|x_1 - x_2| - 1, L - 2).$$
\subsection{Ghost RSA}

We also consider a finite analogue of the ghost RSA process of~\cite{TS06}, which arises from the following Poisson point process that generalizes random sequential adsorption.
We consider a process in $\RR^d$ where the centers of candidate spheres of radius $1$ arrive continuously for time $t \in \RR_{\ge 0}$ according to a translationally invariant Poisson point process of density $\eta = 1$ per unit time (in other words, we expect to see one newly arrived sphere per unit volume and time).

Unfortunately, this process does not create a packing as spheres can overlap. Thus, we must \textit{thin out} the candidate spheres so that the remainder is a packing.

\begin{defn}
For $\kappa \in [0, 1]$, the \textit{$\kappa$-packing} of $\RR^d$ is achieved by running the above described Poisson process, in which we only retain candidate sphere at position $r$ and time $t$ if no other candidate sphere was within a unit distance of $r$ in the time interval $[(1-\kappa)t, t]$. In particular, $\kappa = 0$ corresponds to the \textit{random sequential addition} (RSA) process.
\end{defn}

\begin{defn}
The \textit{ghost RSA process} is the $\kappa$-packing process with $\kappa = 1$.
\end{defn}

For a packing of $\RR^d$, the amount of space occupied by spheres can be quantified by the notions of packing density and number density.

\begin{defn}
Given a packing of $\RR^d$ by spheres of radius $\frac12$, the \textit{packing density} $\phi(t)$ is the fraction of space in $\RR^d$ covered by the spheres. The \textit{number density}  $\rho(t) = \phi(t)/V_d$, where $V_d$ is the volume of the unit sphere in $\RR^d$.
\end{defn}

One of the beautiful properties of the ghost RSA process is that, in sharp constrast to most sphere packing problems (including classical RSA), it is possible to compute the packing density of this process in \textit{all dimensions} $d$. This packing density is also relatively close to the best generic lower bound on the densest sphere packings.

\begin{thm}[Torquato-Stillinger~\cite{TS06}]
The ghost RSA process and the associated underlying Poisson point process enjoy the following properties:
\begin{itemize}
    \item The expected number of candidate centers in a volume $\Omega$ region at time $t$ is $\Omega t$.
    \item The probability that a region of volume $\Omega$ is empty of candidate centers is $\exp(-\Omega t)$.
    \item In the infinite time limit $(t \to \infty)$, we have
    $\EE[\phi(\infty)] = \frac{1}{2^d}$ and $\EE[\rho(\infty)] = \frac{2^d}{V_d}.$
\end{itemize}
\end{thm}

\section{The maximum gap size in the hard-core model}\label{s:maxgap}

In this section, we prove~\cref{t:maxgap}.
Consider the classical 1D hard-core model, where we randomly place rods of length $2$ onto an interval of length $L$  until we no longer can. Recall that $G(L,r)$ denotes the number of gaps of length at least $r$ at saturation.  We seek a threshold $r = r(L)$ such that as $L\to\infty$, we are likely to find gaps smaller than $r$ and unlikely to find any gap of size much greater than $r$. Towards this goal, we prove~\cref{t:maxgap}.

We first consider fixing $r$, and show that as $L \to \infty$, $\EE[G(L,r)]$ converges to a linear function $c_r(L+2)$, by studying a recurrence relation $\EE[G(L,r)]$ satisfies. Since $G(L,r)$ is weakly decreasing with respect to $r$, so must be $c_r$ as a function of $r$. By quantifying the rate of convergence, we obtain a lower bound on $c_{2-\delta}$ by a linear function of $\delta$. This implies that $\EE[G[L,r]]$ changes from $o(1)$ to $\Omega(1)$ as $2-r$ does. We obtain the desired concentration around this expected value by via the second moment method.

Since $G(L,r)=0$ for $r\geq 2$, we henceforth suppose $r \in (0, 2)$. $\EE[G(L,r)]$ satisfies the following recurrence relation:

\begin{obs}
\label{o:expectation-recurrence}
For every $0<r<2$, the expectation
$\EE[G(L,r)]$ satisfies the following integral recurrence relation:
$$\EE[G(L,r)] = 
\begin{cases}
\frac{2}{L-2} \int_{ 0}^{L-2} \EE[G(x,r)]\,dx & L > 2\\
1 & r \le L < 2 \\
0 & L < r
\end{cases}.
$$
\end{obs}

We will first show that as $L \to \infty$, $\EE[G(L,r)] \to c(L+2)$ for some $c = c_r > 0$.  To do this, we consider some helpful auxiliary functions.
\begin{defn}
Fix $r \in (0, 2)$ and let $f_r(L) = \EE[G(L,r)]$. Define 
\[
 g_r(L) := \frac{f_r(L)}{L+2}.
\]
\end{defn}

\begin{obs}\label{o:grec}
For all $r\in(0,2)$, $c\in\RR$, and $L>2$, we have $$
g_r(L)=\frac{2}{L^2-4}\int_{0}^{L-2}g_r(x)(x+2)\,dx.
$$
\end{obs}
\begin{proof}
The relation follows from a direct calculation:
\begin{align*}
    g_r(L)&=\frac{1}{L+2}\cdot \frac{2}{L-2}\int_{x=0}^{L-2}f_r(x)\,dx=\frac{2}{L^2-4}\int_{0}^{L-2}g_r(x)(x+2)\,dx.
\end{align*}
\end{proof}

We can check by hand that for $L \in [0, 4]$ that $g_r(L) \in [0, 1]$. This continues to hold for larger $L$ by studying the above recurrence in~\cref{o:grec}, yielding the following:
\begin{obs}\label{o:fbounds}
For all $r \in (0, 2)$ and $L\geq 0$, we have $0\leq f_r(L)\leq L+2$.
\end{obs}
We now prove that $g_r$ converges to some $c_r$ as $L \to \infty$ by controlling the derivative $g_r'(L)$.

\begin{lemma}\label{l:gprimeiter}
For every $m \in \NN$, there exists some $N_m \in \NN$ such that for all $r \in (0,2)$, whenever $L > 2m + 1$ and $g_r$ is differentiable at $L$, we have $|g_r'(L)| \le N_m/L^m.$
\end{lemma}
\begin{proof}
For all $L$ at which $g_r$ is differentiable, by~\cref{o:grec}, we have
\begin{align*}
    g_r'(L)&=\frac{-4L}{(L^2-4)^2}\cdot \frac{L^2-4}{2}\cdot g_r(L)+\frac{2L}{L^2-4}\cdot g_r(L-2)\\
    &=\frac{2L}{L^2-4} (g_r(L-2)-g_r(L)).
\end{align*}
By~\cref{o:fbounds}, $0 \le g_r(L) \le 1$ for all $L \ge 0$. Therefore, for all $L \ge 3$, we have
$$|g_r'(L)| \le \frac{2L}{L^2-4}\cdot 2\leq \frac{N_1}{L}$$
for some $N_1 >0$. Substituting this inequality again into the above expression, we see that for all $L \ge 5$, we have
$$|g_r'(L)| \leq \frac{2L}{L^2-4}\cdot \frac{2N_1}{L} \le \frac{N_2}{L^2}$$
for some $N_2 > 0$.
Iterating gives the desired result.
\end{proof}

\begin{lemma}\label{l:limit}
For every $r \in (0, 2)$, there exists $c_r\geq 0$ such that
$\lim_{L \rightarrow \infty} g_r(L) = c_r$. Moreover, the convergence is uniform in $r$.
\end{lemma}
\begin{proof}
Since $g_r$ is differentiable almost everywhere, for all $L_2>L_1>5$, we have 
$$|g_r(L_2) - g_r(L_1)| = \left|\int_{L_1}^{L_2} g_r'(x)\,dx \right| \le N_2 \left(\frac{1}{L_1} - \frac{1}{L_2} \right).$$
Consequently, for all $\eps > 0$, there exists $L_\eps > 0$ such that for all $L_1, L_2 > L_\eps$ we have $|g_r(L_2) - g_r(L_1)| < \eps$. This implies that  $\lim_{L\rightarrow \infty} g_r(L) = c_r$ for some $c_r\geq 0$. Since $N_2$ does not depend on $r$, the convergence is uniform in $r$.
\end{proof}

\cref{l:gprimeiter} actually implies a stronger result. Since 
\[
g_r(L) = c_r- \int_{ L}^{\infty} g_r'(x)\,dx = c_r + o(L^{-m})
\]
for arbitrary large $m \in \NN,$ we have the following:
\begin{cor}\label{c:fastconverge}
For all $m \in \NN$, we have
$$f_r(L) = \left(c_r + o(L^{-m})\right)(L+2),$$
where the convergence is uniform in $r$.
\end{cor}

Given the limiting coefficient $c_r$, we wish to understand its magnitude as a function of $r$. To this end, we define the following auxiliary functions.

\begin{defn}
For every $r \in (0, 2]$, define function $h_r(L)$ with domain $(2,\infty)$ as follows:
$$h_r(L) = \begin{cases}
0 & 2 < L < 2 + r \\
\frac{2}{L-2} & 2 + r \le L \le 4 \\
\frac{2}{L-2}\left(1+\int_{2}^{L-2}h_r(x)\,dx\right) & L > 4.
\end{cases}$$
\end{defn}

For all $L>4$, $f_r(L)$ is continuously differentiable with respect to $r$ on $(0,2)$, with $\frac{\partial f_r(L)}{\partial r}=-h_r(L)$. Moreover, the left derivative of $f_r(L)$ with respect to $r$ at 2 equals $-h_r$. An analysis similar to~\cref{l:gprimeiter} gives the following asymptotic expression for $h_r$.
\begin{lemma}
For every  $r \in (0, 2]$, there exists $\lambda_r>0$ such that $\lim_{L\to\infty}h_r(L)/(L+2)=\lambda_r$,
where the convergence is uniform in $r$.
Moreover, if $r_1<r_2$, then $\lambda_{r_1}\geq \lambda_{r_2}$.
\end{lemma}

Since  $h_r(L)=-\frac{\partial f_r(L)}{\partial r}$ for all $L>4$, one might imagine that $-\lambda_r$ is the derivative of $c_r$ with respect to $r$. We make this notion precise below.
\begin{lemma}\label{l:cderiv}
$\frac{\partial c_r}{\partial r}=-\lambda_r$ for all $r\in(0,2)$, and the left derivative of $c_r$ at 2 equals $-\lambda_2$.
\end{lemma}
\begin{proof}
Consider $\left\{\frac{f_r(L)}{L+2}:L>4\right\}$ and $\left\{-\frac{h_r(L)}{L+2}:L>4\right\}$ as families of functions of $r$. Recall that for all $L>4$, $f_r(L)$ is continuously differentiable with respect to $r$ on $(0,2)$ and that we have the following properties:
\begin{itemize}
    \item $\lim_{L\rightarrow\infty}\frac{f_r(L)}{L+2}=c_r$;
    \item $\frac{\partial}{\partial r}\left(\frac{f_r(L)}{L+2}\right)=-\frac{h_r(L)}{L+2}$;
    \item $-\frac{h_r(L)}{L+2}$ converges to $-\lambda_r$ uniformly in $r$.
\end{itemize}
Applying the differentiable limit theorem, we find that $c_r$ is differentiable on $(0,2)$, with $\frac{\partial c_r}{\partial r}=-\lambda_r$. Extending to $r=2$ gives that the left derivative of $c_r$ at 2 equals $-\lambda_2$.
\end{proof}

The above will be enough understand at what $r$ the expected number of gaps of size $r$ drops from $\Omega(1)$ to $o(1).$

\begin{cor}\label{l:ceps}
Consider function $\gamma(L)$ with image in $(0, 2]$. Then,
\[
\lambda_2(2-\gamma(L))\leq c_{\gamma(L)}.
\]
Moreover, if the image of $\gamma$ lies in $[1,2]$, then we also have
\[
c_{\gamma(L)}\leq \lambda_1(2-\gamma(L)).
\]
Therefore, $c_{\gamma(L)}(L+2)=o(1)$ if and only if $\gamma(L)=2-o(1/L)$.
\end{cor}
\begin{proof}
Since $c_2 = 0$,~\cref{l:cderiv} implies that for all $r\in(0,2)$, we have $c_r=\int_{r}^{2}\lambda_{s}ds$. The desired result then follows by noting that $\lambda_r$ is decreasing with respect to $r$. 
\end{proof}

By applying Markov's inequality, we obtain one side of the threshold from the above. To show the other direction, we will need a second moment result.

\begin{defn}
For $r \in (0, 2)$ and $L>0$, let $V_r(L):=\text{Var}[G(L,r)]$ be the variance of the number of gaps of length at least $r$ on a interval of length $L$ in a uniformly random saturated configuration arising from the hard-core process.
\end{defn}

Observe that $V_r(L)\leq \EE[G(L,r)^2]$. For $0 \le L < r$, $V_r(L) = 0$ and for $r \le L < 2$, $\EE[G(L,r )^2] = 1$. For $L > 2$, we have the following recursive inequality.
\begin{obs}\label{o:variance-recurrence}
For all $r \in (0, 2)$, the following recurrence inequality holds for $V_r(L)$ when $L > 2$:
$$V_r(L) \le \EE[G(L ,r)^2] \le \frac{4}{L-2}\int_{ 0}^{L-2} \EE\left[G(x,r)^2\right]\,dx.$$
\end{obs}
\begin{proof}
For $L>2$, let $X$ denote the left endpoint of the first placed rod, so that $X$ follows the uniform distribution on $[0,L-2]$. Let $G_1(L,r)$ denote the number of gaps of size at least $r$ on the left of the first placed rod, and $G_2(L,r)$ denote the number those on the right. Notice that $G(L,r)=G_1(L,r)+G_2(L,r)$. Further, $G_1(L, r), G_2(L, r)$ are conditionally independent given $X$. We therefore have
\begin{align*}
    \EE[G(L,r)^2]&=\EE\left[(G_1(L,r)+G_2(L,r))^2\right] \\
    &= \frac{1}{L-2} \int_{ 0}^{L-2} \EE\left[(G_1(L,r)+G_2(L,r))^2 \mid X=x\right]\,dx\\
    &\leq\frac{2}{L-2} \int_{ 0}^{L-2} \left(\EE\left[G_1(L,r)^2 \mid X=x\right]+ \EE\left[G_2(L,r)^2 \mid X=x\right]\right)\,dx\\
    &=\frac{2}{L-2}\int_{ 0}^{L-2} \left(\EE\left[G(x,r)^2\right]+ \EE\left[G(L-2-x,r)^2\right]\right)\,dx\\
    &=\frac{4}{L-2}\int_{ 0}^{L-2} \EE\left[G(x,r)^2\right]\,dx.
\end{align*}
\end{proof}

Repeating the previous argument on $\EE[G(L,r)^2]$ instead of $\EE[G(L,r)]$, we obtain the following result analogous to \cref{l:ceps}:

\begin{lemma}\label{l:vub}
There exists constant $\mu_1>0$ such that for any function $\gamma:(0,\infty)\rightarrow[1,2]$, we have
\[
\EE[G(L,r)^2]\leq \mu_1(2-\gamma(L))(L+2)).
\]
\end{lemma}

\begin{proof}[Proof of~\cref{t:maxgap}]
Recall that $f_r(L)=\EE[G(L,r)]$, denoting the expected number of gaps of size at least $r$ in a random saturated configuration.
By~\cref{c:fastconverge}, we have that 
 $$f_r(L) = c_r(L+2) + o(1).$$
If $r = 2 - o(1/L)$, by~\cref{l:ceps} $c_{r}(L+2) = o(1)$. Consequently, the expected number of gaps of size at least $r$ is $o(1)$ and thus by Markov's inequality, with high probability there are no gaps of size at least $r.$

Next fix some $a>0$. By~\cref{l:vub}, there exists $\eps > 0$ (e.g. $\eps=2a\mu_1$) such that $V_{2 - a/L}(L) \le \EE[G(L,2 - a/L)^2]\le \eps$ for $L$ sufficiently large. Meanwhile, by~\cref{l:ceps}, we have $c_{2 - a/L}(L+2) \ge \lambda_2a(L+2)/L\geq \lambda_2a$, and thus for $L$ sufficiently large, we have $f_{2 - a/L}(L) \ge \delta$ for some $\delta=\delta(a)>0$. By the second moment method, we see that
$$\P(G(L,2 - a/L) > 0) \ge \frac{\EE[X]^2}{\EE[X^2]} = \frac{f_r(L)^2}{V_{r}(L) + f_r(L)^2} = \frac{\delta^2}{\eps + \delta^2} \ge \delta'.$$
Consequently, with positive probability, we have a gap of size at least $2 - a/L$ on an interval of length $L$.

Next, take $\ell = \ell(L) = o(L)$ such that $\ell = \omega(1)$. Consider some saturated configuration of a length $L$ interval. For each $i \in [\ell]$, we must have some rod with left endpoint $x_i \in [iL/\ell - 2, i L/\ell]$.  
Every possible choice of $(x_1, \ldots, x_{\lfloor L/\ell - 1\rfloor})$ yields a division of the interval into $\ell$ subintervals of length in between $L/\ell-4$ and $L/\ell$, whose numbers of gaps of size at least $r$ are mutually independent.

By the previous argument, we can choose constants $a, \delta = \delta(a) > 0$ such that for sufficiently large $L$, with probability at least $\delta$, a saturated interval of length greater than $\frac{L}{2\ell}$ has a gap of size at least $2 - a \ell/L$. Consequently, the probability that our length $L$ interval has a gap of size at least $2 - a\ell/L$ is at least 
$$1 - (1 - \delta)^{\ell} = 1-o(1),$$
since $\ell \to \infty$ with $L$.
By choosing $\ell = \log L$, we find that with high probability, an interval of length $L$ has a gap of size at least $2 - a\log L/L$. For all $\eps > 0$ and $L$ sufficiently large,  $2 - a\log L/L > 2 - 1/L^{1-\eps}$, giving the desired result.
\end{proof}

\color{black}

\section{Maximum gaps in the  one-dimensional ghost hard-core model}\label{s:ghostgap}

Here, we prove~\cref{t:ghostgap}, giving a threshold for the maximum gap size in the infinite time limit of the ghost hard-core model.

Consider some iteration of the ghost hard-core model on a length $L$ interval. Imagine that we pause at some $t \in \NN$, and for some choice of $\ell \ll L$, consider the collection of gaps of length $\Theta(\ell)$. For each such gap of size $\Theta(\ell)$, we will attempt to compute the probability this gap is retained as $t \to \infty$ by an inductive argument on the lengths of the segments (including ghosts) that are adjacent to the gap. We make this idea more precise below.

\begin{defn}For $\ell>0$ and $k_1,k_2\geq 0$ such that $k_1+k_2\leq \ell$, let $\P^{(\ell)}(k_1, k_2)$ be the probability that a gap of the form $[x,x+\ell]$, in which $[x,x+k_1]$ and $[x+\ell-k_2,x+\ell]$ are already occupied by ghosts, is eventually retained.
\end{defn}

Note that if $\ell - (k_1 + k_2) \le 2$, then $\P^{(\ell)}(k_1, k_2) = 1$.
\begin{obs}
We have the following recurrence for $\P^{(\ell)}(k_1, k_2)$:
\footnotesize 
\begin{align*}
\P^{(\ell)}(k_1, k_2) &= \left(1 - \frac{\ell- k_1 - k_2 + 2}{L} \right) \P^{(\ell)}(k_1, k_2) + \frac{1}{L} \int_{ 0}^2 \P^{(\ell)}(k_1 + k, k_2) dk + \frac{1}{L} \int_{ 0}^2 \P^{(\ell)}(k_1 , k_2 + k) dk \\
&= \frac{1}{\ell-k_1-k_2 + 2} \int_{ 0}^2 \P^{(\ell)}(k_1+k, k_2)dk + \frac{1}{\ell-k_1 -k_2 + 2} \int_{ 0}^2 \P^{(\ell)}(k_1 , k_2 + k) dk .
\end{align*}
\normalsize
Since $\P^{(\ell)}(k_1,k_2)$ is symmetric in $k_1, k_2$ and only depends on the sum $k_1+k_2$, we take $s := \ell - (k_1+k_2)$ and define
$$
\P^{(\ell)}(s) :=
\begin{cases}
1 & s \le 2 \\
\frac{2}{s+2}  \int_{ s-2}^s \P^{(\ell)}(x)\,dx & s > 2
\end{cases}.
$$
\end{obs}

We first prove the second half of \cref{t:ghostgap}, namely that for all $\eps>0$, with high probability a gap of size at least $(\log n)^{1-\eps}$ is retained in the ghost hard-core process. To do so, we first give a lower bound on $\P^{(\ell)}(s)$.

\begin{lemma}\label{l:ghostlb}
For all $s \ge 2$, $\P^{(\ell)}(s) \ge s^{-s}$.
\end{lemma}
\begin{proof}
We see that this is true for $s \le 2$ and  can check by hand for $2 \le s \le 4$. We verify by induction for $s \ge 4$.
Note that $x^{-x}$ is convex for $x > 1$. We have 
\begin{align*}
    \frac{s^{-s}}{\int_{s-2}^sx^{-x}dx}\leq \frac{s^{-s}}{2(s-1)^{-s+1}}=\frac{1}{2(s-1)}\cdot \left(1+\frac{1}{s-1}\right)^{-s}\leq \frac{2}{s+2},
\end{align*}
and thus
\begin{align*}
    s^{-s}\leq \frac{2}{s+2}\int_{s-2}^sx^{-x}dx\leq \frac{2}{s+2} \int_{ s-2}^s \P^{(\ell)}(x)\,dx=\P^{(\ell)}(s).
\end{align*}
\end{proof}

We imagine that $L\to\infty$ and choose parameter $\ell = \ell(L)$ such that $\ell = \omega(1)$ but $\ell^{\ell + 1} = o(L)$, so that in particular $(1 - 1/\ell^\ell)^{L/\ell} = o(1)$. We will show that with positive probability for some $t \in \NN$, there are $\Omega(L/\ell)$ disjoint gaps of size at least $\ell$ at time $t$ in the ghost hard-core process, and that as $t \to \infty$, at least one of these is retained if $\ell$ is sufficiently small. To show the first claim, we apply the following theorem about a randomly broken interval:

\begin{thm}[Theorem 2.2~\cite{Hol80}]\label{t:broken-interval}
Suppose an interval of length 1 is broken uniformly at random into $n$ subintervals with lengths $S_1\leq\dots\leq S_n$. Then for every $i \in [n]$ and $r\in \NN$, we have
\[
\EE[S_i^r]=\EE[Y_i^r]\cdot\frac{\Gamma(n)}{\Gamma(n+r)},
\]
where
\[
Y_i=\frac{X_n}{n}+\frac{X_{n-1}}{n-1}+\dots+\frac{X_{n-i+1}}{n-i+1},
\]
and $X_1,\dots,X_n$ are independent exponential random variables with mean 1.
\end{thm}
To simplify the calculations below, we will henceforth assume (without loss of generality) that parameters describing a number of rods placed or some time step of the hard-core process are integers.

\begin{lemma}\label{l:numbiggaps}
Fix arbitrary $c \in (0, 1/e)$ and $\gamma\in(0,1)$. Suppose $\ell = \ell(L)$ is a function of $L$ such that $\ell^{\ell + 1} = o(L)$. For $L$ sufficiently large, 
when $L/\ell$ candidate rods have been placed, with probability at least
\[
\left(1-\frac{1}{(1+\gamma)\ln(1/c)}\right)^2(1-\gamma),
\]
there are at least $cL/\ell$ pairs of adjacent rod centers having distance at least $\ell$.
\end{lemma}
\begin{proof}We rescale and consider placing rods of length $2/L$ on an interval of length $1$. In the recaled setting, we count pairs of adjacent rod centers whose distance is at least $\ell/L$.

Let $M=L/\ell$. At the point where $M$ candidate rods have been placed on the unit interval, let $S_{M-\gamma M}$ denote the $\gamma M$-th largest distance between adjacent pairs of rod centers. By \cref{t:broken-interval}, we have
\begin{align*}
    \EE[S_{M-c M}]&=\frac{1}{M}\sum_{i=c M+1}^{M}\frac{1}{i}\sim \frac{\ln(1/c)}{M},\\
    \EE[S^2_{M-c M}]&=\frac{1}{M(M+1)}\left(\left(\sum_{i=c M+1}^M\frac{1}{i}\right)^2+\sum_{i=c M+1}^M\frac{1}{i^2}\right)\sim \frac{\ln^2(1/c)}{M^2}.
\end{align*}
For $M$ sufficiently large, by Paley-Zygmund, we have
\begin{align*}
    \PP(S_{M-c M}\geq 1/M)\geq \PP\left(S_{M-c M}\geq \frac{\EE[S_{M-c M}]}{(1+\eps)\ln(1/c)}\right)\geq \left(1-\frac{1}{(1+\gamma)\ln(1/c)}\right)^2(1-\gamma).
\end{align*}
\end{proof}

Consider attempting to place a new rod at time $t$. This new rod can reduce the size of at most one existing gap. This implies that given two distinct, disjoint gaps (separated by at least one placed rod), the events that each of these gaps are retained as $t \to \infty$ are independent.
We are now ready to conclude the second half of \cref{t:ghostgap}.

\begin{lemma}\label{l:underlog}
If $\ell = (\log L)^{1- \eps}$ for $\eps \in (0, 1)$, then with high probability, in the infinite time limit of the ghost hard-core model, there is at least one gap of size at least $\ell$ (when packing on an interval of length $L$ for sufficiently large $L)$.

\end{lemma}
\begin{proof}
For any $\delta \in (0, 1)$, we can choose  $c\in(0,1/e)$ and $\gamma\in(0,1)$ via~\cref{l:numbiggaps} such that for $L$ sufficiently large, with probability at least $1-\delta/2$ there are at least $cL/\ell$ gaps of size at least $\ell$ at some point in the ghost hard-core process. Consequently, there are at least $\frac{cL}{2\ell}$ gaps of size within $[\ell,2\ell/c]$. By \cref{l:ghostlb}, the probability that none of these gaps is retained is at most
\[
\left(1-\frac{1}{(2\ell/c)^{2\ell/c}}\right)^{\frac{cL}{2\ell}}\leq \exp[-L/(2\ell/c)^{2\ell/c+1}]=o(1),
\]
because we have that 
\begin{align*}
    \frac{L}{(2\ell/c)^{2\ell/c+1}}&=2^{\log L -((1-\eps)\log\log L-\log c)(c(\log L)^{1-\eps}+1)} \overset{L \to \infty}\longrightarrow \infty
\end{align*}
Hence with probability at least $1-\delta$, there exists a gap of size at least $\ell$ in the infinite time limit. Sending $\delta\to 0$ gives the result.
\end{proof}

To show the first half of \cref{t:ghostgap}, we give an upper bound on $\P^{(\ell)}(s)$.

\begin{lemma}\label{l:pbound}
There exists $C > 0$ such that for all $s \ge 0$, 
$\P^{(\ell)}(s)\leq Cs^{-s/3}$.
\end{lemma}
\begin{proof}
There exists constant $M>2$ such that for all $s>M$, we have 
\begin{align*}
    \frac{s^{-s/3}}{\int_{s-2}^sx^{-x/3}dx}\geq\frac{s^{-s/3}}{2(s-2)^{-s/3+2/3}}=\frac{1}{2(s-2)^{2/3}}\cdot\left(1+\frac{2}{s-2}\right)^{-s/2}\geq\frac{2}{s+2}.
\end{align*}
Take some $C>0$ such that $\P^{(\ell)}(s)\leq Cs^{-s/3}$ for all $0\leq s\leq M$. Then by induction, for all $s>M$, we have
\begin{align*}
    Cs^{-s/3}\geq \frac{2}{s+2}\int_{s-2}^sCx^{-x/3}dx\geq \frac{2}{s+2} \int_{ s-2}^s \P^{(\ell)}(x)\,dx=\P^{(\ell)}(s). 
\end{align*}
\qedhere
\end{proof}

\begin{lemma}\label{l:nolog}
For sufficiently large $L$, with high probability, the largest gap that remains as $t \to \infty$ in the ghost hard-core process on an interval of length $L$ has size less than $\log L$.
\end{lemma}
\begin{proof}
Let $\ell=\log L$. At any time, there are at most $L/\ell$ distinct gaps of length at least $\ell$. Applying~\cref{l:pbound} and a union bound, the probability that at least one of them is retained is at most
$$ \frac{L}{\ell}\cdot C\ell^{-\ell/3}=\frac{CL}{\log L^{\log L/3+1}}=o(1).$$
\end{proof}

\cref{t:ghostgap} then follows by combining~\cref{l:underlog,l:nolog}.

\section{Further directions}
\label{s:discuss}
The higher dimensional analogues of RSA and parking are of particular importance. A primary motivating question is trying to understand the maximum density of a \textit{sphere packing}, a maximum collection of congruent radius one spheres in $\RR^d$ that do not overlap.
Determining the densest packings in arbitrary dimensions is one of the most longstanding open problems in discrete geometry, recently resolved in $\RR^8$ and $\RR^{24}$ in the breakthrough works of~\cite{VIA17,CKM17}. The only other dimensions in which optimal sphere packings are known are dimensions $1, 2,$ and $3$. One can derive lower bounds on an optimal sphere packing by studying packing procedures, such as the random sequential addition process of hard spheres in $\RR^d.$ Further, RSA in more than one dimension is in and of itself a process of much physical interest, as in~\cite{PWN91,SBY07,MAT74,ITOH11,IU83,EVA93,BES70,BES82}. 

\begin{figure}
    \centering
    \begin{subfigure}[b]{0.48\textwidth}
    \centering
    \includegraphics[width=0.95\linewidth]{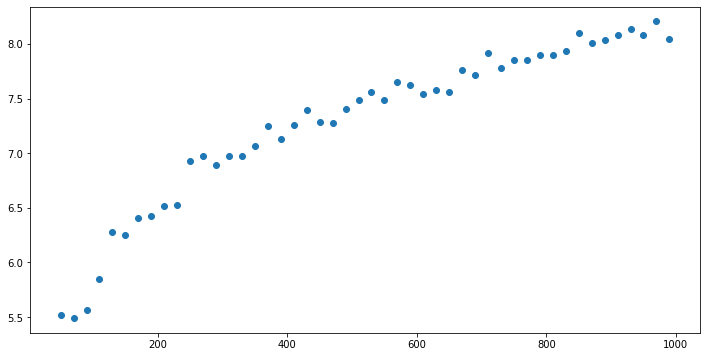}
    \caption{}
    \end{subfigure}
    \begin{subfigure}[b]{0.48\textwidth}
    \centering
    \includegraphics[width=0.95\linewidth]{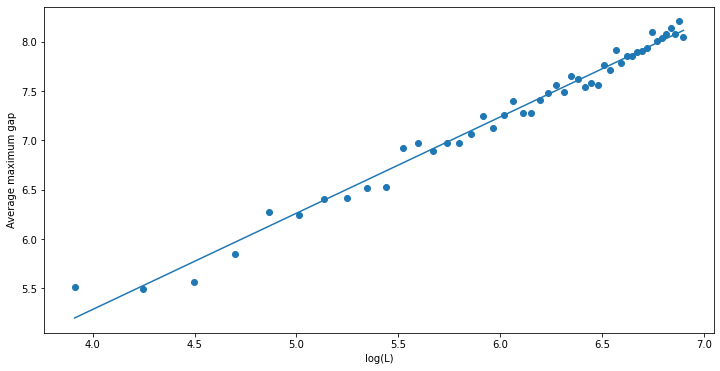}
    \caption{}
    \end{subfigure}    
    \caption{We simulate running the 1D ghost hard-core process to completion, packing rods of width $2$ onto an interval of length $L$ for $L$ ranging from $50$ to $1000$. We estimate the maximum gap size empirically by running 200 trials at each length $L$.  In (a), we plot the average maximum gap vs. $L$; in (b), we plot the average maximum gap vs. $\log(L)$ and include a best fit line that highlights the roughly logarithmic nature of the maximum gap, even at small scale.}
    \label{fig:ghosts}
\end{figure}

While our methods for establishing the extreme values of gaps do not generalize to more than one dimension,~\cref{t:maxgap} offers a tantalizing glimpse into the existence of relatively large gaps in saturated hard-core packings and perhaps studying sphere packing densities. Further the extremely fast convergence of the maximum gap size to a roughly logarithmic scale~\cref{fig:ghosts} provides some evidence about the utility of relatively small scale simulations.

Our work leaves several questions open; perhaps the most fundamental open problem is extending the results in this work to higher dimensions.

\begin{qn}
Given a random packing that results from packing spheres of radius $1$ in $\RR^d$ via the $d$-dimensional ghost RSA process, how much more dense (on average) is the saturated packing that results from adding spheres to this existing packing via the traditional RSA process (i.e., ignoring the ghost constraint)?
\end{qn}

\begin{figure}[t!]
    \centering
    \begin{subfigure}[b]{0.32\textwidth}
    \centering
    \includegraphics[width=0.95\linewidth]{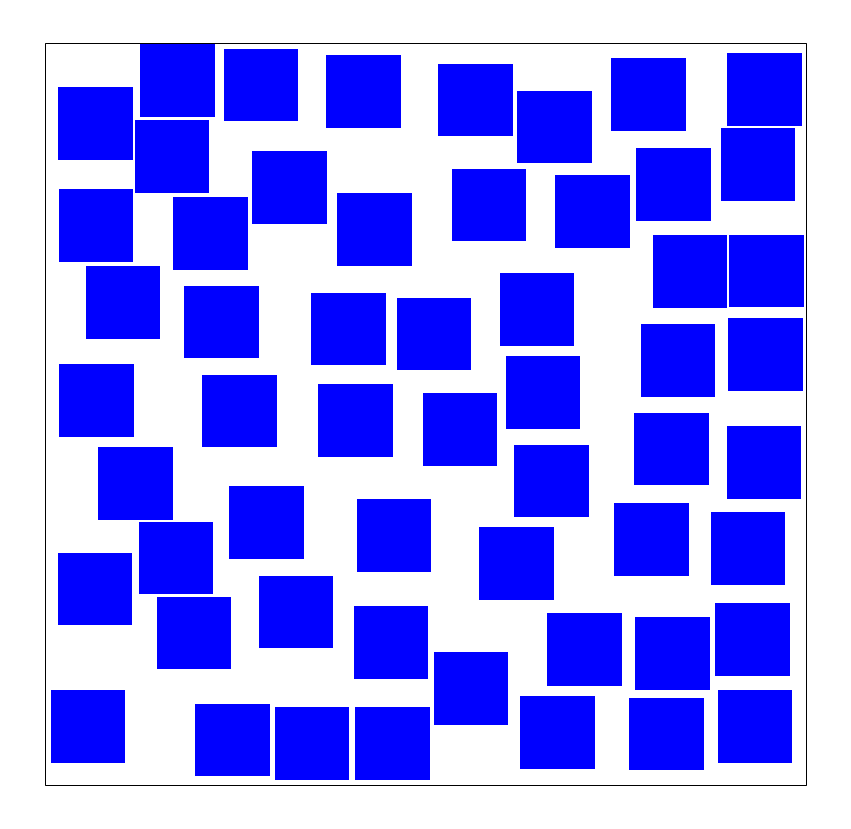}
    \caption{}
    \end{subfigure}
    \begin{subfigure}[b]{0.32\textwidth}
    \centering
    \includegraphics[width=0.95\linewidth]{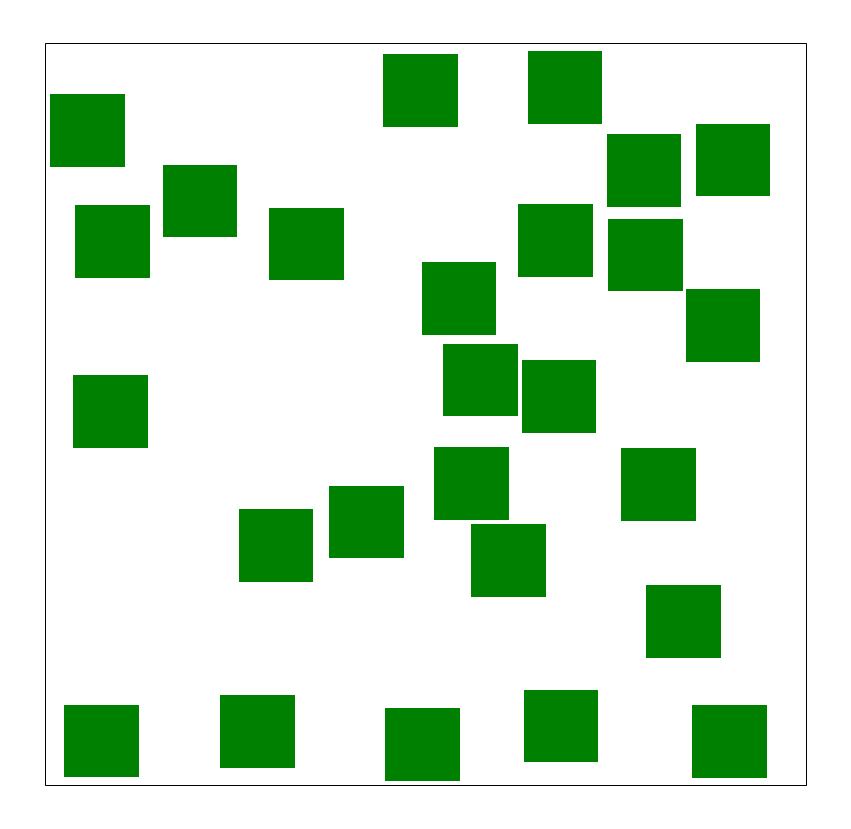}
    \caption{}
    \end{subfigure}
    \begin{subfigure}[b]{0.32\textwidth}
    \centering
    \includegraphics[width=0.95\linewidth]{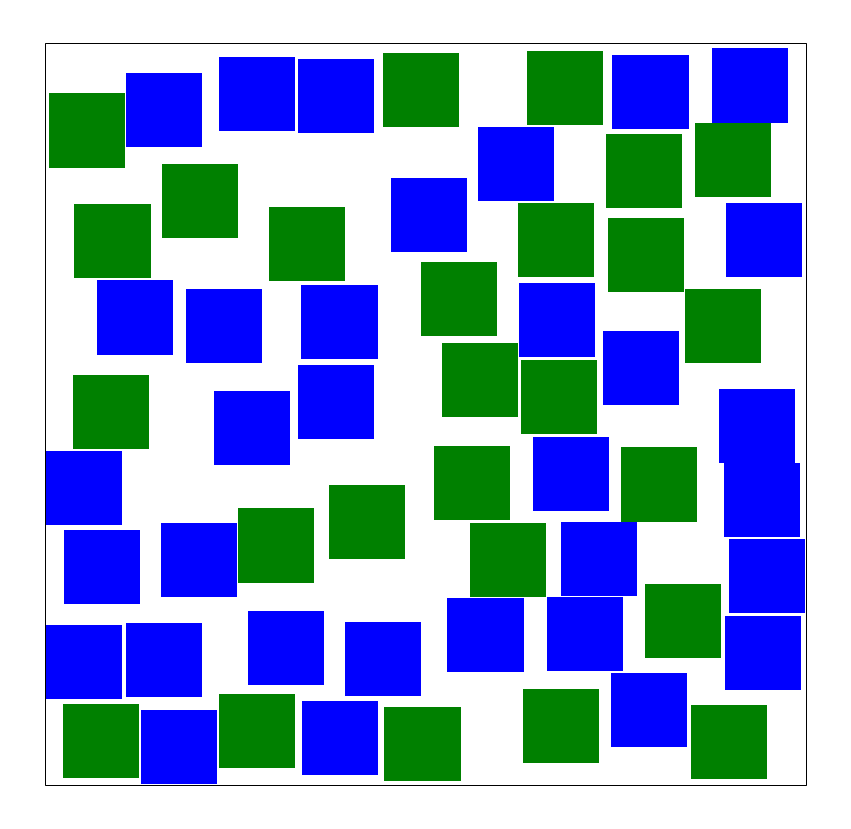}
    \caption{}
    \end{subfigure}
    \caption{Pictured are three ``infinite time'' instantiations of packing axis-aligned squares of side length $2$ inside a larger square of length $L = 20$. In (a), we simulate the classical 2D hard-core model; this packing has 54 squares. In (b), we add the additional \textit{ghost constraint} that squares cannot overlap with any previous candidate square; this packing has 23 squares. In (c), we begin with a random ghost packing generated via the same process as (b); after reaching the infinite time limit,  we extend it randomly to a classically saturated packing by removing the ghost constraint. This example placed 27 squares during the ghost process, and 49 total squares.}
    \label{fig:2d}
\end{figure}

It is also natural to wonder about other generalizations. For example, the following question concerning packing width $2$ axis-aligned squares into an $L \times L$ square is a natural first extension (see~\cref{fig:2d} for a sample simulation).

\begin{qn}
Given a random packing of width $2$ axis-aligned squares into an $L \times L$ square via the hard-core process, what is the expected maximum width of an axis-aligned square that could be added to  this packing, without violating the hard-core constraint?
\end{qn}

Unlike in one dimension, where the notion of a maximum gap is unambiguous, one can define a ``gap'' in many ways in higher dimensions. Focusing on the largest spheres that can fit into the negative space as our notion of gaps, the following question about gaps in the 2D ghost  process is very natural. Two-dimensional packings have been extensively studied; physical theories such as the Asakura-Oosawa depletion interaction seek to explain the ``effective interaction'' between hard-sphere particles and highlight the rich behavior of particles in higher dimensions~\cite{ASA54}. Most of the work in more than one dimension, however, remains at the heuristic level, lacking rigorous mathematical results about distributions of gaps or sizes (or even asymptotic packing densities and pair correlations)~\cite{KRA06}. For example, the following question is open.

\begin{qn}
 Given a random packing of radius 1 spheres into an $L \times L$ square via the ghost hard-core process, what is the expected maximum radius of a sphere that could be added to this packing, without violating the hard-core constraint?
\end{qn}

\bibliographystyle{plainurl}
\bibliography{ghost.bib}

\appendix

\section{Basic properties of the 1D ghost hard-core model}\label{s:basicghost}
Since the \textit{ghost hard-core model} (of~\cref{d:ghost}) is not as well understood as the classical hard-core model, below, we record some basic properties of this ghost model in one dimension below. Some of them bear striking resemblance to the classical results, while other properties of this ghost hard-core process are very different.

In the ghost hard-core process, for each $t \in \ZZ_{\ge 1}$ we attempt to place a length $2$ rod on an interval of length $L$ subject to the \textit{ghost hard-core} constraint that our rod cannot overlap any previous \textit{candidate} rod. Let $\P_s(t)$ be the probability we are successful at time $t$.

\begin{obs}
For $t \in \ZZ_{\ge 1}$, 
\begin{align*}
\P_s(t) &=\frac{L-4}{L} \left(\frac{L-4}{L}\right)^{t-1} + 2 \cdot \frac{1}{L} \int_{1}^2 \left(\frac{L-x-2}{L}\right)^{t-1} dx + 2 \cdot \frac{1}{L} \cdot 0 \\
&= \frac{2 \left((L-3)^t-(L-4)^t\right)
}{t L^t}+\left(\frac{L-4}{L}\right)^t.
\end{align*}
\end{obs}

Recall that $\widetilde{N}(L)$ denotes the total number of rods placed in the ghost model in the infinite time limit (on an interval of length $L$.

\begin{obs}
We have that
$$\EE[\widetilde{N}(L)] = \sum_{t = 1}^{\infty} \P_s(t) = \sum_{t = 1}^{\infty}\left[ \frac{2 \left((L-3)^t-(L-4)^t\right)
}{t L^t}+\left(\frac{L-4}{L}\right)^t\right] = \frac{L}{4} + \ln\left(\frac{16}{9}\right)-1.$$
Thus, the expected density of rods is 
$$\frac{2}{L} \cdot \left( \frac{L}{4} + \ln\left(\frac{16}{9}\right) -1\right) \approx \frac12 - \frac{0.425}{L}.$$
\end{obs}

\begin{rem}
The above is in contrast to the slightly different process of placing length $2$ rods on a circle of length $L$ subject via the one dimensional ghost hard-core model. On a circle of length $L$, the probability of success at time $t \in \ZZ_{\ge 1}$ is simply $\left(\frac{L-4}{L}\right)^{t-1}$. Consequently, the expected density of rods on this circle, is \textit{exactly} $\frac12$ rather than exhibiting the slight offset from $\frac12$ observed in the case of the line.

Naively, one might expect to be able to relate the packing densities on the line segment and circular segment. In the classical one-dimensional hard-core model, we can, after placing the first rod on a circular segment of length $L$ imagine ``cutting'' the rod in half and unrolling it to a line of length $L$ of which a length $1$ piece on the right and left are both occupied. The expected number of rods that can be placed is then just $1$ more than $\EE[N(L-2)]$, the expected number of rods that can be placed on a line segment of length $L-2$ subject to the usual hard-core constraint.

However, the ``unwinding'' argument above does not yield such a correspondence in the ghost hard-core process on a circle; after unwinding, ``ghost'' rods that are placed on top of the first rod can still forbid new space, making the unwound segment exhibit different physical behavior than a length $L-2$ segment.
\end{rem}

\subsection{Occupancy probability}
We now compute the occupancy probability distribution in the infinite time limit of the one-dimensional ghost hard-core process.

\begin{defn}
Let $\widetilde\pi( x,t,L)$ be the probability that point $x \in [0, L]$ is covered by a placed rod for the first time at time $t \in \ZZ_{\ge 1}$, and let $\widetilde \pi(x,L) = \sum_{t = 1}^{\infty} \widetilde\pi(x,t,L)$ be the probability that point $x$ is covered in the infinite time limit.
\end{defn}

\begin{obs}
Suppose that $0 \le x \le \frac{L}{2}$ and $L \ge 10$. Then, $$ \widetilde\pi(x,L) = \sum_{t = 1}^{\infty} \widetilde\pi( x,t,L) = 
\begin{cases}
\frac12 & x \ge 3 \\
\frac{x-1}{4} + \ln\left(\frac{4}{1+x}\right) & 2 \le x < 3 \\
\frac{x-1}{4} + \ln\left(\frac{4}{3}\right) & 1 \le x < 2 \\ 
 \ln\left(\frac{3+x}{3}\right) & 0 \le x < 1.
\end{cases}
$$
\end{obs}
\begin{proof}
This follows by a direct calculation:
\begin{align*}
\widetilde\pi(x,t,L) 
&= \begin{cases}
\frac{2}{L} \left(\frac{L-4}{L}\right)^{t-1} & x \ge 3 \\
\frac{x-1}{L}\left(\frac{L-4}{L}\right)^{t-1} + \frac{1}{L^t} \int_{x-1}^2 (L-y-2)^{t-1} dy & 2 \le x < 3 \\
\frac{x-1}{L}\left(\frac{L-4}{L}\right)^{t-1} + \frac{1}{L^t} \int_{1}^2 (L-y-2)^{t-1} dy & 1 \le x < 2 \\
\frac{1}{L^t} \int_{1}^{x+1} (L-y-2)^{t-1} dy & 0 \le x < 1 \\
\end{cases}    \\
&= 
\begin{cases}
\frac{2}{L} \left(\frac{L-4}{L}\right)^{t-1} & x \ge 3\\
\frac{x-1}{L}\left(\frac{L-4}{L}\right)^{t-1} + \frac{(L-x-1)^t - (L-4)^t}{t L^t} & 2 \le x < 3 \\
\frac{x-1}{L}\left(\frac{L-4}{L}\right)^{t-1} + \frac{(L-3)^t - (L-4)^t}{t L^t}  & 1 \le x < 2 \\
 \frac{(L-3)^t - (L-3-x)^t}{t L^t}  & 0 \le x < 1.
\end{cases}    
\end{align*}
Thus, we sum the above to obtain the occupancy probabilities in the infinite time limit.
$$ \widetilde\pi(x,L) = \sum_{t = 1}^{\infty} \widetilde\pi( x,t,L) = 
\begin{cases}
\frac12 & x \ge 3  \\
\frac{x-1}{4} + \ln\left(\frac{4}{1+x}\right) & 2 \le x < 3 \\
\frac{x-1}{4} + \ln\left(\frac{4}{3}\right) & 1 \le x < 2 \\ 
 \ln\left(\frac{3+x}{3}\right) & 0 \le x < 1.
\end{cases}
$$
\end{proof}

\begin{rem}
The above holds symmetrically for $x > L/2$ by replacing $x$ with $L - x$ in the above expressions.
We plot $\widetilde\pi(x,L)$ for small values of $x$ in Figure~\ref{fig:occup}, noting the boundary effect and the lack of $L$ dependence (provided $L \ge 10$).
\begin{figure}
    \centering
    \includegraphics[width=0.4\textwidth]{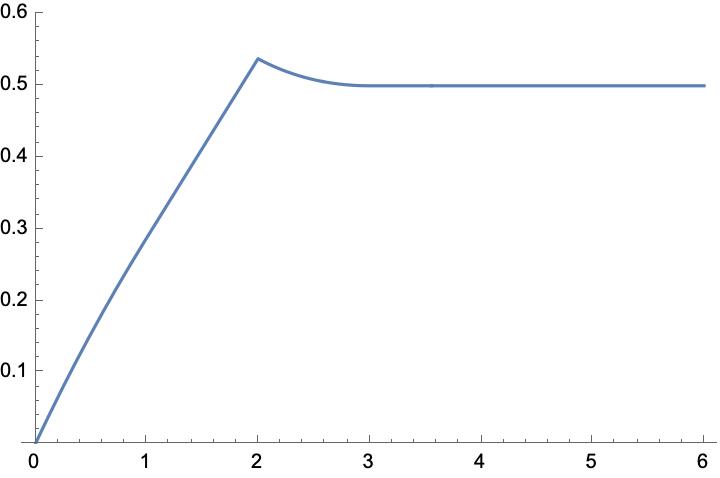}
    \caption{We plot $\widetilde\pi(x)$ for $x \in [0, 6]$ above on an interval of length $L = 20$. $\pi(x)$ is maximized at $x = 2$, where it takes on value $\frac14 + \ln(4/3)$}
    \label{fig:occup}
\end{figure}
\end{rem}

\subsection{Pair correlation}

Let's now consider the pair correlation function on a circular segment of length $L$.
Recall that given a circular segment of length $L$,  $\widetilde \pi(x_1, x_2; L) $ is the probability that  points $x_1$ and point $x_2$ are both covered by a rod when we pack a circle of length $L$ with length 2 rods uniformly at random via the ghost hard-core process. Since $\widetilde \pi(x_1,x_2 ; L)$ only depends on $|x_1-x_2|$, for all $0\leq x\leq L/2
$, we can let $ \widetilde\pi_2(x, L)$ denote $\widetilde \pi(x_1,x_2,L)$ for all $|x_1-x_2|=x$.

\begin{prop}
Let $x = |x_2 - x_1|$. Then we have that for $L \ge 20$ and $x \le L/2$
$$\widetilde \pi_2(x; L) = \widetilde  \pi_2(x) = \frac14 \begin{cases}
1 & x \ge 4 \\
x - \frac{x^2}{8} - 1 & 3 \le x \le 4\\
\frac{19}{8} - \frac12 x - 3 \ln\frac{4}{x+1} & 2 \le x \le 3\\
\frac{x^2}{8}-x+\frac{23}{8}- \ln
  \left(\frac{64}{27}\right) & 1 \le x \le 2 \\
2+\ln(27)-3\ln(x+3)   & 0 < x \le 1.
\end{cases}$$
\end{prop}
%\nitya{I am here!}
\begin{proof}
We imagine a circular segment of length $x = |x_2 - x_1|$ being placed uniformly on the length $L$ circle and compute the probability that in some saturated packing, both endpoints of this segment are within a rod. First consider the left endpoint $x_1$. Since $x_1$ is a uniformly random point on the circle, the probability that $x_1$ is covered by a rod is $\frac12$. We will thus condition on this $\frac12$ event.

By rotational invariance, if $x_1$ is covered, we can assume that $x_1$ was covered at time $t = 1$. 
We can then imagine cutting at point $x_1$ to give a segment of length $L-2$ with a protrusion on the left of length $\sigma$ and on the right of length $2-\sigma$, for $\sigma \sim \text{Unif}([0, 1])$ representing the parts of the cut rod that covered $x_1$. 

For fixed $\sigma,$ let $\widetilde\pi(x, t \mid \sigma)$ be the probability of covering $x$ first at time $t$. Below, we assume $t \ge 2$. We take $x = |x_2 - x_1|,$ assuming that $x_2$ is to the right of $x_1$ and that $x < L/2$. For $t = 1$, 
$$\widetilde\pi(x, 1 \mid \sigma) = 
\begin{cases}
1 & x \le \sigma \\
0 & \text{else}.
\end{cases}
$$
We first suppose that $\sigma \ge 1$. Then we observe that for $t \ge 2$, 
\begin{align*}
\widetilde\pi(x, t \mid \sigma) 
&= 
\begin{cases}
\frac{2}{L} \left(\frac{L - 4}{L}\right)^{t-2} & 2 + \sigma \le x \\
\frac{x - \sigma}{L} \left(\frac{L - 4}{L}\right)^{t-2} & \sigma < x \le 2 + \sigma \\
0 & x \le \sigma.\\
\end{cases}   
\end{align*}
Then the probability of $x$ being covered in the limit for given $\sigma \ge 1$ is
\begin{align*}
\widetilde\pi(x \mid \sigma) &= \sum_{t = 1}^{\infty} \widetilde\pi(x, t \mid \sigma) =
\begin{cases}
\frac12 & 2 + \sigma \le x \\
\frac{x - \sigma}{4} & \sigma < x \le 2 + \sigma \\
1 & x \le \sigma \\
\end{cases}   
\end{align*}

Next we suppose that $\sigma < 1$. We then observe that for $t \ge 2$, we have that 
\begin{align*}
\widetilde\pi(x, t \mid \sigma) &=
\begin{cases}
\frac{2}{L} \left(\frac{L - 4}{L}\right)^{t-2} & 3 \le x \\
\frac{x-1}{L}\left( \frac{L-4}{L} \right)^{t-2} + \int_{x-1}^2 \frac{1}{L} \left(\frac{L-y-2}{L}\right)^{t-2} dy  & 2 + \sigma \le x < 3 \\
\frac{x-1}{L}\left( \frac{L-4}{L} \right)^{t-2} + \int_{\sigma+1}^2 \frac{1}{L} \left(\frac{L-y-2}{L}\right)^{t-2} dy  & 1 \le x < 2 + \sigma \\
\int_{\sigma+1}^{x+1} \frac{1}{L} \left(\frac{L-y-2}{L}\right)^{t-2} dy  & \sigma < x \le 1 \\
0 & x \le \sigma 
\end{cases}    \\
&= 
\begin{cases}
\frac{2}{L} \left(\frac{L - 4}{L}\right)^{t-2} & 3 \le x \\
\frac{x-1}{L} \left(\frac{L - 4}{L}\right)^{t-2} + \frac{(L-x-1)^{t-1} - (L-4)^{t-1}}{(t-1)L^{t-1}} & 2 + \sigma \le x < 3 \\
\frac{x-1}{L} \left(\frac{L - 4}{L}\right)^{t-2} + \frac{(L-\sigma-3)^{t-1} - (L-4)^{t-1}}{(t-1)L^{t-1}} & 1 \le x < 2 + \sigma \\
\frac{(L-\sigma-3)^{t-1} - (L-x-3)^{t-1}}{(t-1)L^{t-1}} & \sigma < x \le 1 \\
0 & x \le \sigma.
\end{cases}   
\end{align*}
The probability of $x$ being covered in the limit for given $\sigma < 1$ is given by
\begin{align*}
\widetilde\pi(x \mid \sigma) &= \sum_{t = 1}^{\infty} \widetilde\pi(x, t \mid \sigma) =
\begin{cases}
\frac12 & 3 \le x \\
\frac{x-1}{4} + \ln \frac{4}{x+1} & 2 + \sigma \le x < 3 \\
\frac{x-1}{4} + \ln \frac{4}{3+\sigma}  & 1 \le x < 2 + \sigma \\
\ln \frac{3+x}{3+\sigma} & \sigma < x \le 1 \\
1 & x \le \sigma.
\end{cases}  
\end{align*}

Thus, taking the minimal distance $|x_1 - x_2| < \frac{L}{2}$ in the regime $L \ge 10$,
\begin{align*}
\widetilde \pi(x_1, x_2, L) &= \widetilde\pi_2(|x_2 - x_1|; L) \\
&:= \widetilde\pi_2(x) \\
&= \frac14 \int_{0}^1 \widetilde\pi(x \mid \sigma) d\sigma + \frac14 \int_1^2 \widetilde\pi(x \mid \sigma) d \sigma  \\
&= \frac14 \begin{cases}
1 & x \ge 4 \\
x - \frac{x^2}{8} - 1 & 3 \le x \le 4\\
\frac{19}{8} - \frac12 x - 3 \ln\frac{4}{x+1} & 2 \le x \le 3\\
\frac{x^2}{8}-x+\frac{23}{8}- \ln
  \left(\frac{64}{27}\right) & 1 \le x \le 2 \\
2+\ln(27)-3\ln(x+3)   & 0 < x \le 1.
\end{cases}
\end{align*}
\end{proof}

\begin{figure}[t]
    \centering
\includegraphics[width=0.4\linewidth]{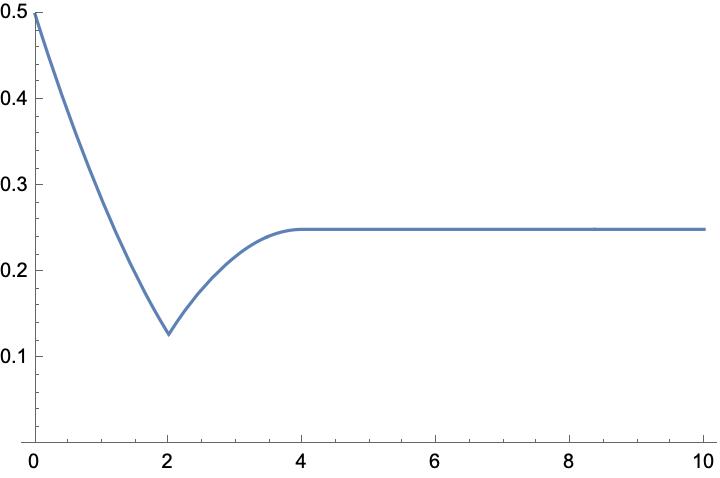}
    \caption{Above, we plot $\widetilde\pi_2(x)$ as a function of $x = |x_2 - x_1|$. The minimum occurs at $x = 2$ with associated occupancy probability $\frac{11}{32} - \frac34 \ln(4/3) \approx 0.128$.}
    \label{fig:pairplot}
\end{figure}

We plot this piecewise function for $x \in [0, 10]$ in Figure~\ref{fig:pairplot}.

\end{document}